\documentclass[a4paper,11pt]{amsart}
\usepackage{amssymb}
\usepackage[english]{babel}
\usepackage[latin1]{inputenc}
\usepackage{graphicx}

\theoremstyle{plain}
\newcommand{\ep}{\epsilon}
\newcommand{\R}{\mathbb R}

\begin{document}
\title{Properties of option prices in models with jumps}
\author[Erik Ekström and Johan Tysk]{Erik Ekström$^1$ and Johan Tysk$^{2,3}$}
\email{ekstrom@maths.manchester.ac.uk, Johan.Tysk@math.uu.se}
\thanks{$^1$ School of Mathematics, The University of Manchester, Sackville Street,
Manchester M60 1QD, UK}
\thanks{$^2$ Department of Mathematics, Box 480, SE-751 06 Uppsala, Sweden}
\thanks{$^3$ Partially supported by the Swedish Research Council}


\newtheorem{theorem}{Theorem}[section]
\newtheorem{lemma}[theorem]{Lemma}
\newtheorem{corollary}[theorem]{Corollary}
\newtheorem{proposition}[theorem]{Proposition}
\newtheorem{definition}[theorem]{Definition}
\newtheorem{hypothesis}[theorem]{Hypothesis}

\newenvironment{example}[1][Example]{\begin{trivlist}
\item[\hskip \labelsep {\bf Example}]}{\end{trivlist}}
\newenvironment{remark}[1][Remark]{\begin{trivlist}
\item[\hskip \labelsep {\bf Remark}]}{\end{trivlist}}

\begin{abstract}
We study convexity and monotonicity properties of option prices in a model with jumps using
the fact that these prices satisfy certain parabolic integro-differential equations.
Conditions are provided under which preservation of
convexity holds, i.e. under which the value, calculated under a chosen martingale measure,
of an option with a convex contract function is convex as a function of the underlying stock price.
The preservation of convexity is then used to derive monotonicity
properties of the option value with respect to the different parameters of the model,
such as the volatility, the jump size and the jump intensity.
\end{abstract}

\maketitle

\section{Introduction}
The implications of parameter mis-specifications for option pricing in continuous diffusion models
have been studied rather extensively in the literature,
compare for example \cite{BGW}, \cite {E}, \cite{EJT},
\cite{EJS}, \cite{H}, \cite{JT1} and the references therein. A general result of these papers
is that in the Black-Scholes model, with the
volatilities of the stocks being possibly time- and level-dependent,
a European option price is monotonically increasing in the volatility if and only if the option
price, at each fixed time prior to maturity, is convex as a function of the price of the
underlying assets. We refer to a model
as being  ``convexity preserving'' if, for any convex contract
function, the corresponding option price is convex as a function of the underlying asset
at all times prior to maturity. Before proceeding, we should note that there are several examples
of models in finance that are not convexity preserving. For example, general stochastic
volatility models, compare \cite{BGW} and \cite{EJS},
and also several models in higher dimensions (even in the case of time- and level-dependent
volatility), compare \cite{EJT}.

In this paper we investigate properties of option prices in a one-dimensional model with possible 
jumps
(for a nice introduction to financial models with jumps we refer to the book \cite{CT}).
When introducing jumps into the model, completeness of the model is in general lost,
so there is not a unique equivalent martingale measure that can be used for arbitrage free
pricing of options. There are at least two different natural questions along the lines of the above
mentioned monotonicity results. First, given a model specified under the physical measure,
which of two given martingale measures gives the higher option price?
Second, what are the implications of a possible mis-specification of models when
using a {\em fixed} martingale measure for pricing options?
The first issue is dealt with in \cite{HH}. In the current paper we focus on the
second question. Thus,
instead of considering different possible choices of martingale measures, we
take the standpoint of having one equivalent martingale measure as given, and we will
investigate properties of the option value calculated under this fixed measure.
This corresponds to an agent modeling the stock
price process, not under the physical measure, but rather directly under a martingale measure,
and with a possible uncertainty about the parameters.

The literature about option price orderings in models with possible jumps is, at least
to the best of our knowledge, not as extensive as the corresponding literature for 
models without jumps.
As mentioned above, Henderson and Hobson \cite{HH} prove
monotonicity results for the option value with respect to different
martingale measures. Their results can also be interpreted
as a monotonicity result in the intensity of the underlying Poisson process for a
{\em fixed} measure, compare Remark 6.2 in that paper.
For the case of American options, see also Pham \cite{P2}.

To remain in a Markovian setting, all parameters of the model are assumed to be
merely time- and level-dependent in this article. Thus the only source of randomness in the
parameters is through their dependence on the current stock price.
In contrast to the case without jumps, it is easy to construct examples of one-dimensional
Markovian models with jumps where convexity is not preserved; such an example is given in
Section~\ref{convexity}.
Below we provide a sufficient condition for a model with jumps to
be convexity preserving. We also show, by analogy with the diffusion case,
that convexity of the value function
implies certain monotonicity results in the different parameters of the
model. More precisely, we derive conditions under which the option value of a convex claim
is increasing in (i) the diffusion coefficient measuring the continuous fluctuations of the
diffusion, (ii) the jump intensity of the underlying Poisson process and (iii) the possible
jump sizes. In the proofs of the preservation of convexity and the monotonicity results
we use a characterization of the option value as the unique viscosity solution of a certain
parabolic integro-differential equation. The methods of proofs are adapted from techniques
used to study preservation of convexity of solutions to parabolic partial differential
equations, see \cite{JT2}.

The present paper is organized as follows. In Section~\ref{def}
we introduce the financial model in the case of finite jump intensity.
In Section~\ref{regularity} a regularity result for the value function of the option
is provided. In Sections~\ref{convexity} and \ref{monotonicity}
we formulate our main results for finite intensity models. 
We give a sufficient condition under which the
model is convexity preserving, and we apply this to establish monotonicity results in
the volatility, the jump size and the jump intensity.
Finally, in Section \ref{activity} we generalize our results 
to models with infinite intensity of jumps.

\section{The financial model}
\label{def}

We consider a financial market with a finite time horizon $T$. The market consists of
a bank account with deterministic interest rate and a risky asset with a positive price process
$X(t)$. For our purposes there is no restriction to assume that the price of the risky asset
is quoted in terms of the non-risky asset, i.e. that the bank account serves as a numeraire. 

The model described below is for the case of finite intensity of jumps; 
this will be generalized in Section 
\ref{activity} to models with possibly infinitely many jumps in finite time intervals.

Let $(\Omega, \mathcal F, P)$ be a probability space with filtration
$\mathbb F=(\mathcal F_t)_{t\in[0,T]}$ with $\mathcal F_T=\mathcal F$ and
satisfying the ``usual conditions'', i.e. $\mathbb F$ is
right-continuous and $\mathcal F_0$ contains all $P$-negligible events in $\mathcal F$.
Let $W$ be a Brownian motion and let $v$ be a Poisson random measure on
$[0,T]\times[0,1]$ with intensity measure
\[q(dt,dz)=\lambda(t)\,dt\,dz\]
for some deterministic function $\lambda(t)\geq 0$.
Define the compensated jump martingale measure $\tilde v$ by $\tilde v(dt,dz)=(v-q)(dt,dz)$.
Let the risky asset be modeled by a stochastic process $X(t)$ satisfying
the stochastic differential equation
\begin{equation}
\label{X}
dX=\beta(X(t-),t)\,dW +\int_0^1\phi(X(t-),t,z)\,\tilde v(dt,dz).
\end{equation}
The interpretation that should be given to the model is that
\begin{itemize}
\item
$\beta$ (or rather $\frac{\vert\beta\vert}{x}$)
represents the volatility of the Brownian part of the stock price;
note that $\beta$ is possibly time- and level-dependent.
\item
the jump intensity of the stock price is $\lambda$; associated to each jump there is a
label $z$ with the interpretation that a jump at time $t$ with label $z$ is of size
$\phi(X(t-),t,z)$.
\end{itemize}
The following are the main assumptions used in Section \ref{def} to Section \ref{monotonicity}
in this paper.
\begin{itemize}
\item[(M1)]
The diffusion coefficient $\beta:\R^+\times[0,T]\to\R$ and the jump intensity $\lambda:[0,T]\to\R$
are both continuous,
and the jump size $\phi:\R^+\times[0,T]\times[0,1]\to\R$ is measurable and for each fixed $z\in[0,1]$,
the function $(x,t)\mapsto\phi(x,t,z)$ is continuous.
\end{itemize}
Moreover, there exist constants $C>0$ and $\gamma>-1$ with
\begin{itemize}
\item[(M2)] $\beta^2(x,t)+\phi^2(x,t,z)\leq Cx^2$
\item[(M3)]
$\vert\beta(x,t)-\beta(y,t)\vert+\vert\phi(x,t,z)-\phi(y,t,z)\vert\leq C\vert x-y\vert$
\item[(M4)]
$\phi(x,t,z)>\gamma x$
\end{itemize}
for all $x$, $t$ and $z$.
Under these assumptions there exists a unique solution to the stochastic differential
equation (\ref{X}) for any starting point $x>0$, compare for example \cite{GS}.
This solution satisfies
\[P\big(X(t)\leq 0\mbox{ for some }t\in[0,T]\big)=0,\]
i.e. $X(t)$ remains strictly positive at all times with probability 1.

For a continuous contract function $g$, the value at time $t$ of a European option
that at time $T$ pays the amount $g(X(T))$ is $u(X(t),t)$, where
\[u(x,t)=E_{x,t}g(X(T)).\]
Here the indices indicate that $X(t)=x$.
In Lemma~3.1 in \cite{P}, estimates of the second moment of $X(T)$ are given.
In our setting with a bounded intensity of jumps, estimates on higher
moments than two can be obtained in the same way as in the well-known case of
diffusion processes using Gronwall's lemma. Thus the value function $u$ is well-defined
for contract functions satisfying
\begin{itemize}
\item[(M5)]
$g\in C_{pol}(\R^+)$,
\end{itemize}
where $C_{pol}(\R^+)$ is defined in Definition~\ref{spaces} below.

\begin{remark}
As remarked in the introduction, there is in general no unique risk-neutral measure that
could be used for arbitrage free pricing of options in the above jump-diffusion model.
In this paper we do not deal with the issue of choosing an appropriate measure for pricing,
but we rather assume that the measure has already been chosen. Moreover,
we specify our model, not under a physical measure, but directly
under this chosen martingale measure. Thus
there is no need of changing measures when pricing options.
\end{remark}

The pricing function $u$ can under appropriate conditions
be characterized as the unique viscosity solution of a certain parabolic
integro-differential equation. In fact, after a standard change $t\to T-t$ of the direction of time,
$u$ satisfies
\[u_t=\mathcal L u\]
with initial condition
\[u(x,0)=g(x).\]
Here $\mathcal L$ is the elliptic integro-differential operator
\begin{equation}
\label{evalundgren}
\mathcal L u:=a u_{xx}+ \mathcal Bu,
\end{equation}
where
\[a(x,t):=\frac{\beta^2(x,t)}{2}\]
and
\[\mathcal Bu=\lambda(t)\int_0^1 \Big(u(x+\phi(x,t,z),t)-u(x,t)-
\phi(x,t,z)u_x(x,t)\Big)\,dz.\]
We will see below that, under some conditions, a viscosity solution also is a classical
solution. Therefore we do not formally introduce the concept of viscosity solutions here, but
we rather refer the reader to the definition in \cite{P}.

\section{Regularity of the value function}
\label{regularity}
In this section we provide some regularity results for the value function $u$.
A priori, viscosity solutions are merely continuous, but by arguing similarly as in Section~5.2
in \cite{P}, higher order regularity can be obtained (see also Theorem 4, page 296 in \cite{GS}).
To do this we first transform the equation into an equation with coefficients satisfying some
standard assumptions in the theory of partial differential equations.
We would like to point out that the results in \cite{P} easily extend to the case
of a possibly time-dependent jump-intensity $\lambda$.

We begin with introducing a few definitions and assumptions that are used below.

\begin{definition}
\label{spaces}
(i)
For a set $E\subset\R$ we denote by $C^p(E)$ the set of functions $f:E\to\R$ such that the derivatives
$\partial_x^k f$ with $k\leq p$ exist in the interior of $E$ and have continuous extensions to $E$.

(ii)
For a set $E\subset \R\times [0,T]$
we denote by $C^{p,q}(E)$ the set of functions $f:E\to\R$ such that the derivatives
$\partial_x^k\partial_t^l f$ with $k+2l\leq p$ and $l\leq q$ exist in the interior of $E$
and have continuous extensions to $E$.

(iii)
For a set $E\subset\R\times [0,T]$
we denote by $C_{pol}(E)$ the set of functions of at most polynomial growth in $x$.
More explicitly,
\[C_{pol}(E)=\bigcup_{C>0,m>0}\big\{f:E\to\R:\vert f(x,t)\vert\leq
C(1+\vert x\vert^m)\mbox{ for }(x,t)\in E\big\}.\]

(iv)
For $E\subset\R\times[0,T]$ and $\alpha\in(0,1)$
we denote by $C_\alpha(E)$ the set of locally H\"older($\alpha$) functions,
i.e.
\[C_\alpha(E):=\{f:E\to\R:\sup_{p,q\in K}
\frac{\vert f(p)-f(q)\vert}{d(p,q)^\alpha}<\infty\mbox{ for each compact $K\subset E$}\}.\]
Here $d$ is the parabolic distance $d\big((x_1,t_1),(x_2,t_2)\big)=(\vert x_1-x_2\vert^2+
\vert t_1-t_2\vert)^{1/2}$.

(v)
For $E\subset \R\times [0,T]$, the spaces $C^{p,q}_{pol}(E)$ and $C^{p,q}_\alpha(E)$ are
the spaces of functions $f\in C^{p,q}(E)$
for which all the derivatives $\partial_x^k\partial_t^l f$ with $k+2l\leq p$ and $l\leq q$
belong to $C_{pol}(E)$ and $C_\alpha(E)$, respectively.

(vi)
For a set $E\subset\R$, the spaces $C_{pol}(E)$, $C_{pol}^{p,q}(E)$, $C_\alpha(E)$ and
$C^{p,q}_\alpha(E)$ are defined similarly.
\end{definition}

To prove the regularity result Theorem~\ref{reg} below, we need some additional assumptions
beyond the main assumptions (M1)-(M5). We will say that assumptions (A1)-(A5) hold if
there exist constants $C$, $\gamma>0$ and $\alpha\in(0,1)$ such that
\begin{itemize}
\item[(A1)]
$\gamma x^2\leq\beta^2(x,t)$
for all $t$ and all $x\leq\gamma$;
\end{itemize}
\begin{itemize}
\item[(A2)]
$\lambda\in C_\alpha([0,T])$;
\item[(A3)]
$\beta\in C^{2,1}_\alpha(\R^+\times[0,T])$ with
$\vert \beta_t(x,t)\vert\leq Cx$ and $\vert \beta_{xx}(x,t)\vert\leq C/x$;
\item[(A4)]
$\phi(\cdot,\cdot,z)\in C^{2,1}_\alpha(\R^+\times [0,T])$ (with the H\"older continuity
being uniform in $z$), and
\[\vert\phi_t(x,t,z)\vert\leq Cx\]
and
\[\vert\phi_{xx}(x,t,z)\vert\leq Cx^{-1}\]
for all $(x,t,z)$;
\item[(A5)]
$g$ is Lipschitz continuous, i.e. $\vert g(x_2)-g(x_1)\vert\leq C\vert x_2-x_1\vert$
for all $x_1,x_2\in\R^+$. Moreover, $g\in C^3_{pol}(\R^+)$.
\end{itemize}

\begin{theorem}
\label{reg}
In addition to the main assumptions (M1)-(M5), assume that (A1)-(A5) hold.
Then the value function $u$ is in $C^{4,1}(\R^+\times(0,T))\cap C^{2,1}(\R^+\times [0,T])$.
Moreover, there exist constants $m>0$ and $K>0$ such that
\[\vert u_{xx}(x,t)\vert\leq K(x^{-m}+x^m)\]
for all $(x,t)\in\R^+\times [0,T]$.
\end{theorem}

\begin{proof}
Let $Y(t):=\Psi\big(X(t)\big)$ where $\Psi:\R^+\to\R$ is a smooth
function with $\Psi(x)=-1/x$ for $x\in(0,1]$, $\Phi(x)=x$ for $x\geq 2$ and $\Psi^\prime(x)>0$
for all $x\in\R^+$.
Applying the Ito formula for diffusions with jumps, compare for example Chapter 8 in \cite{CT},
it follows that $Y$ solves the stochastic differential equation
\[dY=\tilde b\big(Y(t-),t\big)\,dt+\tilde\beta\big(Y(t-),t\big)\,dW +
\int_0^1\tilde\phi\big(Y(t-),t,z\big)\,\tilde v(dt,dz)\]
in $\R\times [0,T]$, where
\begin{eqnarray*}
\tilde b(y,t) &=& \frac{1}{2}\Psi^{\prime\prime}\big(\Psi^{-1}(y)\big)\beta^2\big(\Psi^{-1}(y),t\big)\\
&&+
\lambda(t)\int_0^1 \Big(\tilde\phi(y,t,z)-y
-\Psi^\prime\big(\Psi^{-1}(y)\big)\phi\big(\Psi^{-1}(y),t,z\big)\Big)\,dz,
\end{eqnarray*}
\[\tilde\beta(y,t)=\Psi^\prime\big(\Psi^{-1}(y)\big)\beta\big(\Psi^{-1}(y),t\big)\]
and
\[\tilde\phi(y,t,z)=\Psi\Big(\Psi^{-1}(y)+\phi\big(\Psi^{-1}(y),t,z\big)\Big)-y,\]
and $\tilde v$ is the same compensated jump martingale measure as in (\ref{X}).
It is straightforward to check that the coefficients $\tilde b$, $\tilde\beta$ and $\tilde \phi$ together
with the initial condition
\[\tilde g(y):=g(\Psi^{-1}(y))\]
satisfy the conditions (2.2)-(2.5) in \cite{P}. Therefore, it follows from Theorem~3.1 in \cite{P} that
\[v(y,t):=u\big(\Psi^{-1}(y),t\big),\]
in the viscosity sense, solves the equation
\[v_t=\frac{\tilde\beta^2}{2}v_{yy}+\tilde bv_y+\lambda\int_0^1 v(y+\tilde\phi,t)-v-\tilde\phi v_y\, dz\]
in $\R\times (0,T]$, with initial condition $v(y,0)=\tilde g(y)$.
Reasoning as on page 22 in \cite{P}, one finds that $v$ also solves,
again in the viscosity sense, the equation
\begin{equation}
\label{ferlin}
\left\{ \begin{array}{l}
v_t=\frac{\tilde\beta^2}{2}v_{yy}+\big(\tilde b-\lambda\int_0^1\tilde\phi\,dz\big) v_y+h\\
v(y,0)=\tilde g(y),\end{array}\right.
\end{equation}
where
\[h(y,t)=\lambda(t)\int_0^1 \big(v(y+\tilde\phi,t)-v(y,t)\big)\, dz.\]
Moreover, Proposition~3.3 in \cite{P} gives that the function $v(y,t)$ satisfies
\begin{equation}
\label{kalifatides}
\vert v(y_1,t_1)-v(y_2,t_2)\vert\leq C\Big((1+\vert y_1\vert)\vert t_1-t_2\vert^{1/2}+
\vert y_1-y_2\vert\Big)
\end{equation}
for some constant $C$ (Theorem~3.1 and Proposition~3.3 in \cite{P} are stated for a problem
of optimally stopping of a controlled jump diffusion process, but the corresponding proofs also work
in our simpler setting).

Using (\ref{kalifatides}) and the assumptions on $\phi$ it follows that $h$ is in
$C_\alpha(\R\times [0,T])\cap C_{pol}(\R\times [0,T])$.
From Theorem A.20 in \cite{JT2} we thus have the existence of a unique solution
$w$ to equation (\ref{ferlin}) with  $w\in C^{2,1}_{pol}(\R\times [0,T])$ (note that
the proof of that theorem
also works with the current weaker condition on $h$), and from Theorem~A.18 in \cite{JT2}
we find that $w\in C^{2,1}_\alpha(\R\times[0,T])$. Moreover, this function $w$ also satisfies
the inequality (\ref{kalifatides}) (this can be seen by noting that the classical solution
$w$ to (\ref{ferlin}) also is a stochastic solution, and thus Proposition~3.3 in \cite{P}
can be applied also to $w$).

Since $w$ is a classical solution to (\ref{ferlin}), it is also a viscosity solution of
this equation. It then follows from the uniqueness result Theorem~4.1 in \cite{P}
that $v\equiv w$. Consequently, $v\in C^{2,1}_{pol}(\R\times [0,T])\cap C^{2,1}_\alpha(\R\times[0,T])$.
It follows that $h\in C^{2,0}_\alpha(\R\times [0,T])$, so applying Theorem~A.11 in \cite{JT2} yields
$v\in C^{4,1}(\R\times(0,T))$.

Changing back to the original coordinates,
it follows that $u\in C^{4,1}(\R^+\times(0,T))\cap C^{2,1}(\R^+\times[0,T])$
and that there exists a constant $m$ such that $u_{xx}=\mathcal O(x^{-m})$ for $x$ close to 0 and
$u_{xx}=\mathcal O(x^{-m})$ for large $x$, uniformly in $t$.
\end{proof}

\begin{remark}
The reason to use the change of coordinates $y=-1/x$ for small $x$, and not the more standard change
$y=\ln x$, is to be able to use Theorem~A.20 in \cite{JT2}. With the logarithmic coordinate
change, the condition (A4) of that theorem would not be fulfilled.
\end{remark}

\section{Preservation of convexity}
\label{convexity}
In this section we provide a sufficient condition on $\phi$ under which the model is
convexity preserving, see condition (\ref{erlander}) below.
The methods used are adapted from \cite{JT2}, in which the same problem is
studied for parabolic equations.
We begin, however, with an example showing that not all models are convexity preserving.

\begin{example}
{\bf (A model which is not convexity preserving.)}
Let $\phi:\R^+\to\R$ be a non-negative Lipschitz function satisfying
$\phi(x)=0$ for $x\notin(1/2,3/4)$
and $\phi(x)>1$ for $x\in(x_1,x_2)$, where $x_1,x_2$ satisfy $1/2<x_1<x_2<3/4$.
Further, let $g(x)=(1-x)^+$ and the stock price dynamics be given by
\[dX=\phi(X(t-))\, (dN-dt),\]
where $N$ is a Poisson process with intensity 1. Now, since $X(t)$ is a martingale,
and since $g$ is convex, it is easy to check that the option value $u(x,0)$ at time 0 satisfies
\[u(x,0)=E_{x,0}g(X(T))\geq g(x),\]
with strict inequality for $x\in(x_1,x_2)$. Moreover, for $x\notin(1/2,3/4)$ the
inequality reduces to an equality, and thus $u(x,0)$ is not convex in the interval $[1/2,1]$.
\end{example}

\begin{definition}
We say that a model is convexity preserving if, for each convex contract function $g\in C_{pol}(\R^+)$,
the corresponding value function $u(x,t)$ is convex in $x$ for each fixed time $t\in[0,T]$.
\end{definition}

Along the lines of the analysis in \cite{JT2}, we make the following definition.
The differential operator $\mathcal L$ is defined as in (\ref{evalundgren}).

\begin{definition}
In addition to the main assumptions (M1)-(M5) of Section~\ref{def}, also assume that (A1)-(A5) hold.
Then we say that the model is locally convexity preserving (LCP) at a point
$(x_0,t_0)\in C(\R^+\times[0,T])$ if
\begin{equation}
\label{LCP}
\partial_x^2(\mathcal L f )(x_0,t_0)\geq 0
\end{equation}
holds for any convex function $f\in C^4(\R^+)\cap C^2_{pol}(\R^+)$ with $f_{xx}(x_0)=0$.
We simply say that a model is LCP if it is LCP at all points.
\end{definition}

\begin{remark}
Note that the condition (M2) ensures the integral term $\mathcal B f$
in $\mathcal L f$ to be
well-defined for any function $f\in C^1(\R^+)\cap C_{pol}(\R^+)$. Similarly, (M2), (M3) and (A4)
together ensure that the integral term in
$\partial_x^2(\mathcal L f)$ is well-defined for any function $f\in C^3(\R^+)\cap C^2_{pol}(\R^+)$.
\end{remark}

It is intuitively clear that the LCP-condition is a natural
condition to impose for preservation of convexity.
Indeed, if spatial convexity of $u$ is almost lost at some point $(x_0,t_0)$, then the infinitesimal
change of $u$ in the time interval $[t_0,t_0+\Delta t]$ is given by
$\Delta t\,\mathcal L u$, which is spatially convex if the LCP-condition is satisfied.
Below we show that a model which is LCP in fact also is convexity preserving.

\begin{theorem}
\label{main}
Let the assumptions of Theorem~\ref{reg} hold, and assume that the model is LCP.
Then the model is convexity preserving.
\end{theorem}

\begin{proof}
We know from Theorem~\ref{reg} that the value function $u$ is in
$C^{4,1}(\R^+\times[0,T])$ and that there exist constants $K$ and $m$ such that
\begin{equation}
\label{schyman}
\vert u_{xx}(x,t)\vert\leq K(x^m+x^{-m})
\end{equation}
for all $(x,t)\in\R^+\times[0,T]$. Define the function $h:\R^+\to\R$ by $h(x):=x^{m+3}+x^{-m+1}$. Then
\begin{equation}
\label{gxx}
h_{xx}=(m+2)(m+3)x^{m+1}+m(m-1)x^{-m-1}
\end{equation}
and
\begin{eqnarray*}
\partial_x^2(\mathcal L h) &=& a h_{xxxx}+2a_xh_{xxx}+a_{xx}h_{xx}\\
&&
+\lambda\int_0^1 \Big((1+\phi_x)^2h_{xx}(x+\phi)+\phi_{xx}h_x(x+\phi)-\phi h_{xxx}\\
&&\hspace{17mm}-(1+2\phi_x)h_{xx}-\phi_{xx}h_x\Big)\,dz.
\end{eqnarray*}
The assumptions (A3)-(A4) on $\beta$ and on $\phi$ imply the existence of a large positive
constant $M$ such that
\begin{equation}
\label{ohly}
M \partial_x^2 h-\partial_x^2(\mathcal L h)\geq 1
\end{equation}
for all $x$ and $t$. For $\ep>0$, define the function $u^\ep$ by
\[u^\ep(x,t):=u(x,t)+\ep e^{Mt}h(x),\]
and assume, to reach a contradiction, that the set
\[E:=\{ (x,t):u\mbox{ is not convex in the spatial variable at }(x,t)\}\]
is non-empty. Since $u_{xx}$ satisfies (\ref{schyman}), it follows from (\ref{gxx}) that
$E\subset (\rho^{-1},\rho)\times [0,T]$ for some $\rho\in\R^+$.
Thus $E$ is bounded, so $\overline E$ is compact. Therefore the infimum
\[t_0:=\inf\{ t\geq 0:(x,t)\in \overline E\mbox{ for some }x\in\R^+\}\]
is attained, and there exists $x_0\in\R^+$ with $(x_0,t_0)\in\overline E$.
At this point we have by continuity $u_{xx}^\ep=0$, and therefore $t_0>0$ (since
$u_{xx}^\ep(x,0)\geq \ep h_{xx}(x)>0$).
By the definition of $t_0$, $u_{xx}^\ep(x_0,t)\geq 0$ for $t\leq t_0$, so
\[\partial_tu^\ep_{xx}(x_0,t_0)\leq 0.\]
Moreover, since $u^\ep_{xx}(x_0,t_0)=0$ and $u^\ep$ is spatially convex at $t=t_0$,
the LCP-assumption yields
\[\partial_x^2(\mathcal L u^\ep)\geq 0\]
at $(x_0,t_0)$. Thus we find that
\[\partial_x^2\big(\partial_tu^\ep-\mathcal L u^\ep\big)(x_0,t_0)\leq 0.\]
On the other hand, using (\ref{ohly}) and the fact that $u$ solves the equation
$u_t=\mathcal L u$, we have
\[\partial_x^2\big(\partial_tu^\ep-\mathcal L u^\ep\big) =
\ep e^{Mt}\partial_x^2 (M h-\mathcal L h)\geq \ep e^{Mt}>0,\]
so we have reached a contradiction. Therefore the set $E$ is empty, and thus $u^\ep$ is spatially
convex at all times. By letting $\ep$ tend to 0 it follows that also $u$ is spatially convex,
finishing the proof.
\end{proof}

\begin{theorem}
\label{branting}
Let the assumptions of Theorem~\ref{reg} hold. Also assume that
\begin{equation}
\label{palme}
\phi_{xx}(x,t,z)\phi(x,t,z)\geq 0
\end{equation}
for all $x,t,z$. Then the model is LCP, and thus also convexity preserving.
\end{theorem}

\begin{proof}
Let $f\in C^4(\R^+)\cap C^2_{pol}(\R^+)$. Then,
\begin{eqnarray*}
\partial_x^2(\mathcal L f) &=& af_{xxxx}+2a_xf_{xxx}+a_{xx}f_{xx}\\
&&+\lambda\int_0^1 \Big((1+\phi_x)^2f_{xx}(x+\phi)+\phi_{xx}f_x(x+\phi)-\phi f_{xxx}(x)\\
&&\hspace{15mm}-(1+2\phi_x)f_{xx}(x)-\phi_{xx}f_x(x)\Big)\,dz.
\end{eqnarray*}
Now, assuming that $f$ is convex and satisfies $f_{xx}(x_0)=0$ at some point $x_0$,
$f_{xx}$ has a local minimum at $x_0$. Thus $f_{xxx}(x_0)=0$ and $f_{xxxx}(x_0)\geq 0$, so, at
a point $(x_0,t_0)$,
\begin{eqnarray*}
\partial_x^2(\mathcal L f)
&=& af_{xxxx}\\
&& +\lambda\int_0^1 \Big((1+\phi_x)^2f_{xx}(x_0+\phi)+\phi_{xx}f_x(x_0+\phi)
-\phi_{xx}f_x(x_0)\Big)\,dz\\
&\geq& \lambda\int_0^1 \Big(\phi_{xx}f_x(x_0+\phi)-\phi_{xx}f_x(x_0)\Big)\,dz.
\end{eqnarray*}
It follows from the assumption (\ref{palme}) and the convexity of $f$ that
\[\phi_{xx}f_x(x+\phi)-\phi_{xx}f_x(x)\geq 0\]
for all $(x,t,z)$. Consequently $\partial_x^2(\mathcal L f)\geq 0$ at the point
$(x_0,t_0)$, so the model is LCP.
\end{proof}

By approximation we can relax the conditions of Theorem~\ref{branting} as follows.

\begin{theorem}
\label{goranpersson}
Assume the main conditions (M1)-(M5) and that
\begin{equation}
\label{erlander}
\phi \mbox{ is convex (concave)
in $x$ at all points where $\phi(x,t,z)>0$ ($<0$).}
\end{equation}
Then the model is convexity preserving.
\end{theorem}

\begin{proof}
Let $(\beta,\phi, \lambda)$ be a model satisfying (M1)-(M4) and (\ref{erlander}), and let
$X$ and $u$ be the corresponding stock and option prices, respectively.
First choose a contract function $g$ satisfying (A5), and
let $(\beta_n,\phi_n,\lambda_n)$ be a sequence of models satisfying (M1)-(M4) uniformly in $n$, i.e.
(M2)-(M4) hold for all models with the same constants $C$ and $\gamma$. Also assume that
each model $(\beta_n,\phi_n,\lambda_n)$ satisfies the conditions
(A1)-(A4) (not necessarily uniformly in $n$) and (\ref{erlander}) (or equivalently, (\ref{palme})),
and that for each $N>0$ and $t\in[0,T]$ we have
\[\lim_{n\to\infty}\sup_{x\in(0,N]}\vert\beta_n(x,t)-\beta(x,t)\vert=0,\]
\[\lim_{n\to\infty}\sup_{x\in(0,N]}\int_0^1\vert \phi_n(x,t,z)-\phi(x,t,z)\vert^2\,dz=0,\]
and
\[\lim_{n\to\infty}\vert \lambda_n(t)-\lambda(t)\vert=0.\]
Then the option price $u_n$ corresponding to the model $(\beta_n,\phi_n,\lambda_n)$
is spatially convex by Theorem~\ref{main},
and the stock price $X_n(T)$, corresponding to the model $(\beta_n,\phi_n,\lambda_n)$,
converges in $L^2$ to $X(T)$, compare Part II, \S8, Theorem 3 in \cite{GS}
(this theorem is stated for a constant jump-intensity, but it readily extends to our setting).
This implies that $u_n(\cdot,t)$ converges pointwise to $u(\cdot,t)$.
Since the pointwise limit of a sequence of convex functions is convex, we find that
$u(\cdot,t)$ is convex.

Finally, for a convex contract function $g$ not necessarily satisfying the Lipschitz
condition (A5), but merely the weaker condition (M5) allowing polynomial growth,
one may approximate $g$ with a sequence of Lipschitz functions.
It is straightforward to check that the corresponding prices converges to the correct limit,
i.e. the result about preservation of convexity extends to
contract functions g of possibly polynomial growth.
\end{proof}

\section{monotonicity in the model parameters}
\label{monotonicity}

In this section we demonstrate how preservation of convexity can be used to derive monotonicity
properties of the option value with respect to the different parameters of the model.
To do this we consider two different models, i.e. two sets $(\beta,\phi, \lambda)$ and
$(\tilde\beta,\tilde\phi,\tilde\lambda)$ of parameters, and we denote by $u$, $\tilde u$ and
$\mathcal L$, $\tilde {\mathcal L}$ the corresponding option values and integro-differential
operators, respectively.

\begin{theorem}
\label{mon}
Assume that both models satisfy the main assumptions (M1)-(M4) stated in Section~\ref{def},
and that the contract function $g$ satisfies (M5).
Also assume that $\vert\tilde\beta(x,t)\vert\leq \vert\beta(x,t)\vert$ and
$\tilde\lambda(t)\leq\lambda(t)$ for all $x$ and $t$, and that
\begin{equation}
\label{cond}
\frac{\phi(x,t,z)}{\tilde\phi(x,t,z)}\geq 1
\end{equation}
for all $x,t,z$ with $\tilde\phi(x,t,z)\not= 0$.
If either $\phi$ or $\tilde\phi$ satisfies condition (\ref{erlander}),
and if the contract function $g$ is convex, then
\[\tilde u(x,t)\leq u(x,t)\]
for all $x$ and $t$.
\end{theorem}

\begin{remark}
Theorem~\ref{mon} extends a result of Henderson and Hobson \cite{HH}. In Theorem~6.1 of that
paper it is shown that, for any convex contract function $g$, the corresponding
option price $u$ is increasing in the intensity $\lambda$ provided all parameters
of the model are deterministic.

Also note that a consequence of Theorem~\ref{mon} is that if the contract function $g$ is convex, then
the Black-Scholes price (corresponding to $\phi\equiv 0$) gives a lower bound for the
set of possible arbitrage free option prices. This is also proven by Bellamy
and Jeanblanc \cite{BJ}.
\end{remark}

\begin{proof}
We assume the assumptions (A1)-(A4) to hold for both models and that $g$ satisfies (A5);
the general case then follows
by an approximation argument similar to the one in the proof of Theorem~\ref{goranpersson}.
Under these assumptions it follows from Theorem~\ref{reg} that there exist large positive
numbers $m$ and $K$ such that
\begin{equation}
\label{peralbin}
\max\{\vert u(x,t)\vert,\vert\tilde u(x,t)\vert\}\leq K(x^m+x^{-m})
\end{equation}
for all $(x,t)\in\R^+\times[0,T]$.
Let $h=x^{m+1}+x^{-m-1}$, and choose the constant $M$ large so that
\begin{equation}
\label{bardot}
Mh-\mathcal L h\geq 1
\end{equation}
for all $x$ and $t$. Define
\[u^\ep(x,t):=u(x,t)+\ep e^{Mt}h(x),\]
and suppose that the set
\[E:=\{(x,t)\in\R^+\times[0,T]:u^\ep(x,t)<\tilde u(x,t)\}.\]
is non-empty. It follows from (\ref{peralbin}) that there exists $\rho>0$ with
$E\subseteq (\rho^{-1},\rho)\times [0,T]$. Thus $E$ is bounded, so $\overline E$ is compact.
Hence there exists a point $(x_0,t_0)\in\overline E$ where
\[t_0=\inf\{ t:(x,t)\in \overline E\mbox{ for some }x\in(0,\infty)\}.\]
By continuity, $u^\ep(x_0,t_0)=\tilde u(x_0,t_0)$, so
$u^\ep(x,0)-\tilde u(x,0)\geq\ep h(x)>0$ implies that $t_0>0$. It is therefore clear that
\begin{equation}
\label{strindberg}
\partial_t(u^\ep-\tilde u)\leq 0
\end{equation}
at the point $(x_0,t_0)$. On the other hand, at this point we also have
\begin{eqnarray}
\label{lagerlof}
\notag
\partial_t(u^\ep-\tilde u) &=&
\mathcal L u^\ep-\tilde{\mathcal L} \tilde u +\ep e^{Mt_0}(M h-\mathcal L h)\\
\notag
&=&
\frac{\beta^2}{2}u^\ep_{xx}-\frac{\tilde\beta^2}{2}\tilde u_{xx}+\ep e^{Mt_0}(Mh-\mathcal L h)\\
\notag
&& +\int_0^1 \Big(\lambda(t_0)\big(u^\ep(x_0+\phi,t_0)-u^\ep(x_0,t_0)-\phi u^\ep_x(x_0,t_0)\big)\\
\notag
&& \hspace{12mm} - \tilde\lambda(t_0)\big(\tilde u(x_0+\tilde\phi,t_0)-\tilde u(x_0,t_0)-\tilde\phi
\tilde u_x(x_0,t_0)\big)
\Big)\,dz\\
&>&
(\frac{\beta^2}{2}-\frac{\tilde\beta^2}{2})u^\ep_{xx}
+\frac{\tilde\beta^2}{2}(u^\ep_{xx}-\tilde u_{xx})\\
\notag
&& +\int_0^1 \Big(\lambda(t_0)\big(u^\ep(x_0+\phi,t_0)-u^\ep(x_0,t_0)-\phi u^\ep_x(x_0,t_0)\big)\\
\notag
&& \hspace{12mm}- \tilde\lambda(t_0)\big(\tilde u(x_0+\tilde\phi,t_0)-\tilde u(x_0,t_0)-\tilde\phi \tilde u_x(x_0,t_0)\big)\Big)\,dz,
\end{eqnarray}
where we have used the inequality (\ref{bardot}). Assume first that $\phi$ satisfies
(\ref{erlander}). From Theorem~\ref{main} it follows that $u$ is spatially convex, and therefore
also $u_\ep$ is spatially convex. Using $u^\ep=\tilde u$ and $u^\ep_x=\tilde u_x$
at $(x_0,t_0)$ and the condition (\ref{cond}) we find that
\[u^\ep(x_0+\phi,t_0)-u^\ep(x_0,t_0)-\phi u^\ep_x(x_0,t_0)
\geq \tilde u(x_0+\tilde\phi,t_0)-\tilde u(x_0,t_0)-\tilde\phi \tilde u_x(x_0,t_0).\]
Since the expression on the left hand side of this inequality due to convexity
is non-negative, we also have
\begin{eqnarray*}
&&\lambda(t_0)\Big( u^\ep(x_0+\phi,t_0)-u^\ep(x_0,t_0)-\phi u^\ep_x(x_0,t_0)\Big)
\geq \\
&&\hspace{10mm}\tilde\lambda(t_0)\Big(\tilde u(x_0+\tilde\phi,t_0)-\tilde u(x_0,t_0)-
\tilde\phi \tilde u_x(x_0,t_0)\Big),
\end{eqnarray*}
so it follows from $\beta^2\geq \tilde\beta^2$ and $u^\ep_{xx}(x_0,t_0)\geq \tilde u_{xx}(x_0,t_0)$ that
\[\partial_t(u^\ep-\tilde u)>0.\]
This contradicts (\ref{strindberg}). Thus the set $E$ is empty, and so $\tilde u(x,t)\leq u^\ep (x,t)$
for all $x$ and $t$. Now the desired monotonicity result follows by letting $\ep\to 0$.

If instead $\tilde\phi$ satisfies (\ref{erlander}), then
$\tilde u$ is spatially convex, so we can argue similarly as above if the expression in
(\ref{lagerlof}) is replaced by
\[\frac{\beta^2}{2}(u^\ep_{xx}-\tilde u_{xx})+
(\frac{\beta^2}{2}-\frac{\tilde\beta^2}{2})\tilde u_{xx}.\]
\end{proof}

\section{Models with infinite intensity of jumps}
\label{activity}

In this section we extend the results of the previous sections
to the case of models with possibly an infinite jump activity. 
To introduce such a model, let $m$ be a Radon measure on $(0,\infty)$, i.e. a positive
measure which assigns finite 
measure to any compact subset of $(0,\infty)$. Then let $v$ be a Poisson random 
measure on $[0,T]\times (0,\infty)$ with intensity measure $\lambda(t)\,dt\,m(dz)$.
The stock is now modelled by
\begin{equation}
\label{Xg}
dX=\beta(X(t-),t)\,dW +\int_0^\infty\phi(X(t-),t,z)\,\tilde v(dt,dz)
\end{equation}
where $\tilde v(dt,dz)=v(dt,dz)-\lambda(t)\,dt\,m(dz)$ is the compensated jump martingale
measure. This model allows for an infinite intensity of jumps since the measure $m$
can be chosen so that the parameter space $(0,\infty)$ has infinite mass. Note that the parameter
$\lambda =\lambda(t)$ in this setting is not the standard intensity of jumps of $X$.

To incorporate these more general models we replace the assumptions (M1)-(M5) with the 
assumptions (G1)-(G5). 

\begin{itemize}
\item[(G1)]
Assumption (G1) agrees with (M1) but with the $z$-parameter space $[0,1]$ replaced by $(0,\infty)$.
\end{itemize}
Moreover, there exists a constant $C>0$ with
\begin{itemize}
\item[(G2)] $\beta^2(x,t)+\int_0^\infty\phi^2(x,t,z)\,m(dz)\leq Cx^2$
\end{itemize}
and
\begin{itemize}
\item[(G3)]
$(\beta(x,t)-\beta(y,t))^2+\int_0^\infty (\phi(x,t,z)-\phi(y,t,z))^2\,m(dz)
\leq C (x-y)^2$
\end{itemize}
for all $x$ and $t$. For simplicity we assume 
\begin{itemize}
\item[(G4)]
$\phi$ is such that the solution $X$ to (\ref{Xg}) stays strictly positive during the 
time interval $[0,T]$ with probability 1.
\end{itemize}
In addition to these assumptions, we also impose
\begin{itemize}
\item[(G5)]
For any positive integer $p$ there exists a constant $C_p$ such that 
\[\int_0^\infty\vert\phi(x,t,z)\vert^p\,m(dz)\leq C_p(1+x^p)\]
for all $x$ and $t$.
\end{itemize}
As a consequence of (G5), moments of all orders of $X(T)$ exist, compare Part II, §7, Theorem 1 in
\cite{GS} where this is stated under slightly different conditions. Therefore 
the option price
\[u(x,t)=E_{x,t}g(X(T))\]
is finite for any contract function $g$ of at most polynomial growth. Note that 
the integrability condition (G5) is implicitly satisfied in the case of finite 
jump intensity studied in Section~\ref{def}. Indeed, (M2) and the finiteness of the
Lebesgue measure for the parameter space $[0,1]$ of $z$ imply (G5).

To formulate our result in the present setting, we again consider two different models, 
i.e. two sets $(\beta,\phi, \lambda)$ and
$(\tilde\beta,\tilde\phi,\tilde\lambda)$ of parameters,
and we denote by $u$ and $\tilde u$ the corresponding option values.

\begin{theorem}
\label{tham}
Assume that both models satisfy the general assumptions (G1)-(G5)
and that the contract function $g$ is convex and of at most polynomial growth.
\begin{itemize}
\item[i)] Assume that $\phi$ satisfies (\ref{erlander}). Then 
$u$ is convex in $x$ at each fixed time $t\in[0,T]$.
\item[ii)]
Assume that $\vert\tilde\beta(x,t)\vert\leq \vert\beta(x,t)\vert$ and
$\tilde\lambda(t)\leq\lambda(t)$ for all $x$ and $t$, and that
\[\frac{\phi(x,t,z)}{\tilde\phi(x,t,z)}\geq 1\]
for all $x,t,z$ with $\tilde\phi(x,t,z)\not= 0$.
If either $\phi$ or $\tilde\phi$ satisfies condition (\ref{erlander}), then
\[\tilde u(x,t)\leq u(x,t)\]
for all $x$ and $t$.
\end{itemize}
\end{theorem}

\begin{proof} To prove i), note that it follows from (\ref{erlander}) and (G2) that for fixed
$t$ and $z$, the function $\phi(x,t,z)$ satisfies $\lim_{x\to 0}\phi=0$ $m$-almost everywhere
and is thus either convex and non-negative or concave
and non-positive as a function of $x$. 
We define $\phi_n$ for $z\in[1/n,n]$ by
\begin{eqnarray*}
&&\phi_n(x,t,z):=\\
&&\hspace{5mm}\left\{\begin{array}{ll}
\sup\{f(x):\mbox{$f$ convex, $f(\cdot)\leq \phi(\cdot,t,z)$ and 
$f_x^+ \leq n$}\} & \mbox{if $\phi>0$}\\
\inf\{f(x):\mbox{$f$ concave, $f(\cdot)\geq \phi(\cdot,t,z)$ and 
$f_x^+ \geq \frac{1-n}{n}$}\} & \mbox{if $\phi<0$},
\end{array}\right.
\end{eqnarray*}
where $f^+_x$ denotes the right derivative. For $z\notin[1/n,n]$ we define
$\phi(x,t,z)=0$.
Since $m$ is a Radon measure we have $m([1/n,n])<\infty$, so this construction gives a model 
with finite jump intensity. Moreover, $\phi_n$ satisfies the conditions (M1)-(M5) and also 
the condition (\ref{erlander}). Therefore $u_n$, the option price corresponding to $\phi_n$, 
is convex as a function of $x$ by Theorem~\ref{goranpersson}.
Since the sequence $(\phi_n)_{n=1}^\infty$ satisfies 
\[\int_0^\infty\vert\phi_n(x,t,z)-\phi(x,t,z)\vert^2\,m(dz)\to 0\]
as $n$ tends to infinity, Theorem 3 in 
§8 \cite{GS} gives that $u_n$ tends to $u$ pointwise by the same argument
as in the proof of Theorem~\ref{goranpersson}. Thus also $u$ is convex in $x$.

A similar argument yields the monotonicity part of the theorem.
\end{proof}

\begin{remark}
Also for options of American type, the results corresponding to Theorems~\ref{goranpersson},
\ref{mon} and \ref{tham} hold. This can be seen by approximating the American option price with a sequence
of Bermudan option prices, all of which are spatially convex and increasing in the
model parameters. The details work as in \cite{E}, where models without jumps are treated.
\end{remark}

\end{document}